\documentstyle[12pt]{article}
\begin{document}
\centerline{\bf Zerofree region for exponenetial sums}

\vskip 0.1in

\centerline{\bf R. BALASUBRAMANIAN}

\vskip 0.5in

\noindent \S1 We consider the following two closed sets in $C^n$. One is 
the diagonal D given by $ (z, z, z, \cdots z_)$. The other is $A = 
\left\{(z_1,z_2,z_3, \cdots z_n): \right .$ \newline 
$\left . e^{z_1} + e^{z_2} +e^{z_3} 
+ 
\cdots + 
e^{z_n}=0\right\}$. Clearly $D \cap A$ is empty. One can ask 
what is 
the 
distance between them.

\vskip 0.1in

In this connection, Stolarsky [1] proved that the distance $d$ is given 
by
$d^2 = (\log \ n)^2 + O (1)$.  Some simple calculations  will make one 
believe
that the point $(k, 0, 0, \cdots 0)$  with $k = \log \ (n-1) + \pi i$ 
which lies on $A$ is one of the closest point to the diagaonal. 
We prove that this is
indeed the case, atleast for sufficiently large $n$. This gives
$d^2 = |k|^2  (1-1/n)$.

\vskip 0.1in

\noindent \S2. We first make the following observations.

\vskip 0.1in
 
\indent If $P (z_1, z_2, \cdots z_n)$ is any point in $C^n$, then the 
nearest point in $D$ to $P$ is $(z, z, z, \cdots z)$ where 
$z = 
(z_1 + z_2 + z_3 + \cdots z_n)/n$. Thus if $P (z_1, z_2, z_3, \cdots 
z_n)$ is a point of $A$  closest to $D$, so are the points 
$(\bar{z_1},\bar{z_2}, 
\cdots 
\bar{z_n})$, $(z_1-a, z_2-a, z_3-a, \cdots z_n-a)$ and $(z_{\sigma(1)}, 
z_{\sigma(2)}, 
z_{\sigma (3)}, \cdots z_{\sigma(n))}$ for any $a$ in $C$ and for any
permutation $\sigma$ of $(1,2,3, \cdots n)$.

\vskip 0.1in

Thus  one  may assume that 

$$
\begin{array}{rc}
    e^{z_1} + e^{z_2} + e^{z_3} + \cdots e^{z_n} = 0  &  (1) \\
\\
    Im (z_n) \geq 0                      &       (2) \\
\\
   \sum^n_{j=1} z_j = 0                  &       (3) \\
\\
   |z_1| \leq |z_2| \leq \cdots \leq |z_n|  &    (4) \\
\\
   \mbox{Further we can assume that} \\
\\
   Im z_n \leq \pi                         &       (5) \\
\end{array}
$$

\vskip 0.1in

since otherwise the point $(z_1, z_2, \cdots  z_{n-1}, z_{n} - 2 \pi i)$ 
is nearer to $D$

\vskip 0.1in

Now consider the point $Q (b_1, b_2, b_3, \cdots b_n)$ where 
$b_j = -k/n$ for $1 \leq j \leq (n-1)$ and

$$
  b_n =  k-k/n
$$

\vskip 0.1in

\noindent Then $Q$ satisfies the conditions above and is at a distance 
$|k|^2 (1-1/n)$ from $D$.

\vskip 0.1in

\noindent Our proof will be complete if we prove that there is no 
point at a  smaller distance.

\vskip 0.1in

\noindent Thus we assume $(1), (2), (3), (4), (5)$ \ and 

$$
\begin{array}{rc}
  |z_1|^2 + |z_2|^2 + |z_3|^2 + \cdots |z_n|^2 \leq |k|^2 (1-1/n) & (6) 
\\
\\
   \mbox{Define}  \ a_j = z_j - b_j &                         (7) 
\\
\\
  \mbox {our aim is to prove that}  \ a_j = 0  \ \mbox{for all} \ j\\
\\
  \mbox{We note that} \ \sum^n_{j=1} \ a_j = 0 &          
(8) 
\end{array}
$$

\vskip 0.1in

\noindent {\bf Lemma 1}  We have \ $\displaystyle{\sum^n_{j=1}} |a_j|^2 
= 
O 
(\log n |a_n|)$

\vskip 0.1in

\noindent {\bf Proof}  We have, from equations (6) and (7)

\vskip 0.1in

$$\sum^n_{j=1} |a_j + b_j|^2 \leq \sum^n_{j=1}|b_j|^2$$

\vskip 0.1in

Hence $\sum |a_j|^2 \leq -2 Rl \sum (\bar{b}_j a_j) = 2 |\sum b_j a_j)|$

\vskip 0.1in

\noindent Substituting the value of $b_j$ and using (8) we get

\vskip 0.1in

$\sum |a_j|^2 \leq 2 |k| |a_n|$ and hence the lemma.

\vskip 0.1in

\noindent {\bf Lemma 2}  We have  $\displaystyle{\sum^{n-1}_{j=1}} 
(e^{a_j} - 1) = 
(n-1) (e^{a_n} -1)$

\vskip 0.1in

\noindent {\bf Proof}  We have $0 = \displaystyle{\sum^n_{j=0}}  e^{z_j} 
=
\sum e^{(b_j+ a_j)}$

\vskip 0.1in

\noindent Now we  substitute the value of $b_j$ to get 
$\displaystyle{\sum^{n-1}_{j=1} e^{a_j} = (n-1) e^{a_n}}$.  The lemma 
follows.

\vskip 0.1in

\noindent {\bf Lemma 3} We have $\displaystyle{\sum^{n-1}_{j=1}} 
(e^{a_j} 
- a_j -1) = 
(n-1)(e^{a_n} -1) + {a_n}$ 

\vskip 0.1in

\noindent {\bf Proof} This  follows from lemma2 and (8) 

\vskip 0.1in

\noindent {\bf Lemma 4}  We have $\sum |a_j|^2 = 0 (\log^2 n)$
\vskip 0.1in

\noindent {\bf Proof}   Since $|a_j| \leq |z_j|+|b_j|$, this follows 
from equation (6) and definition of  $b_j$.

\vskip 0.1in

\noindent Now define $M = \max |a_j|$, the maximum over $1 \leq j \leq 
n-1$.

\vskip 0.1in

\noindent {\bf Lemma 5} We have, $M \leq  0.75 \log n$ 

\vskip 0.1in

\noindent {\bf Proof}  We have , if $j \leq (n-1)$, then  $2|z_j|^2 = 
\leq
\displaystyle{\sum^n_{j=1}} 
|z_j|^2 \leq (\log n)^2+ \pi^2$  by equation (4) and (6). Hence 
 $|z_j| \leq  0.72 \log n$.

\vskip 0.1in

\noindent Since $|a_j| \leq |z_j|+|b_j|$, the lemma follows

\vskip 0.1in

\noindent {\bf Lemma 6} We have $\displaystyle{\sum^{(n-1)}_{j=1}}
e^{a_j}-1-a_j)=O ( n^ {0.75})$

\vskip 0.1in

\noindent {\bf Proof}  Let $S_r = a_1^r + a_2^r + \cdots 
a^r_{(n-1)}$.

\vskip 0.1in

\noindent Then $S_1 = -a_n$ and for $l \geq 2, |S_l| \ \leq \ S_2  
M^{(l-2)}$ where $M$ is the maximum of $a_j$ which is $0.75 \log \ n$ in 
our case.

\vskip 0.1in

\noindent  Now $e^{a_j} - 1 - a_j = \displaystyle{\sum^{\infty}_{l=2}} 
a_j^l/l!$. 

Summing over $j=1$ to $(n-1)$

$\displaystyle{\sum^{(n-1)}_{j=1}} (e^{a_j}-1-a_j) = 
\displaystyle{\sum^{\infty}_{l=2}} S_l/l!$

\vskip 0.1in

Now using the bound  on $M$ and noting that $S_2$ is 
$O(\log \ n)^2)$ by Lemma 4, the result follows.

\vskip 0.1in

\noindent {\bf Lemma 8}. We have $a_n = O (n^{-0.25})$

\vskip 0.1in

\noindent {\bf Proof}  By Lemmas 3 and 6, we have  
    $(n-1) |e^{a_n}-1| = O (n^{0.75}) + O (|a_n)|$.

\vskip 0.1in

\noindent Since $a_n = O (\log \ n)$ by Lemma 4, we get $e^{a_n}-1$ is 
small. This means, $a_n$ is small, upto a multiple 
of $2 \pi i$. Now the assumption (2) and (5) ensure that the multiple 
is zero.

\vskip 0.1in

\noindent Now Lemma 8, combined with Lemma 1 shows that all $a_j$ are 
small. With this information, we relook at Lemma 3. Now the right side 
of Lemma 3 is $>> n |a_n|$ and the left side is $<< 
\displaystyle{\sum^{n-1}_{j=1}} 
|a_j|^2$ which is $O (\log \ n |a_n|)$ by Lemma 1. Thus  $n|a_n| << 
 \log \ n  |a_n|$. 
This easily leads to a contradiction unless $a_n=0$. Then by Lemma 1, 
all $a_j$ are zero. Hence the result.

\vskip 0.1in

\noindent Reference:

\vskip 0.1in

[1]: Stolarsky; K.B: Zero free regions of exponential sums.  Proc. 
A.M.S.89 (1987) (486-488).
\end{document}